# A BAYESIAN $\chi^2$ TEST FOR GOODNESS-OF-FIT

BY VALEN E. JOHNSON

*University of Michigan*

This article describes an extension of classical $\chi^2$ goodness-of-fit tests to Bayesian model assessment. The extension, which essentially involves evaluating Pearson's goodness-of-fit statistic at a parameter value drawn from its posterior distribution, has the important property that it is asymptotically distributed as a $\chi^2$ random variable on $K - 1$ degrees of freedom, independently of the dimension of the underlying parameter vector. By examining the posterior distribution of this statistic, global goodness-of-fit diagnostics are obtained. Advantages of these diagnostics include ease of interpretation, computational convenience and favorable power properties. The proposed diagnostics can be used to assess the adequacy of a broad class of Bayesian models, essentially requiring only a finite-dimensional parameter vector and conditionally independent observations.

**1. Introduction.** Model assessment presents a challenge to Bayesian statisticians, one that has become an increasingly serious problem as computational advances have made it possible to entertain models of a complexity not considered even a decade ago. Because diagnostic methods have not kept pace with these computational advances, practitioners are often faced with the prospect of interpreting results from a model that has not been adequately validated.

Numerous solutions to this problem have been considered. The most orthodox of these depend on the specification of alternative models and the use of Bayes factors for model selection. This approach is reasonable when both a relatively broad class of models can be specified as alternatives, and when implied Bayes factors can be readily computed. Unfortunately, it often happens in practice that neither requirement is satisfied, making this approach impractical for routine application. Complicating the situation still further is the fact that Bayes factors are not defined when improper priors are used in model specification, although this difficulty may be partially circumvented through the use of intrinsic Bayes factors or related devices [e.g., Berger and Pericchi (1996) and O'Hagan (1995)].

A second strategy for assessing model adequacy centers on the use of posterior-predictive model checks. This approach was initially proposed by Guttman (1967) and Rubin (1984), and was extended to more general discrepancy functions by Gelman, Meng and Stern (1996). [Gelfand (1996) has advocated related techniques based on cross-validatory predictive densities.] The primary advantage of posterior-predictive model assessment is its relative ease of implementation. In many models, the output from numerical algorithms used to generate samples from the posterior distribution can be used to generate observations from the predictive model, which in turn can be used to compute $p$-values for the discrepancy function of interest. Posterior-predictive model assessment also facilitates case-diagnostics, which, in many circumstances, are more telling in examining model fit than are global goodness-of-fit statistics. However, such approaches also have an important disadvantage. As Bayarri and Berger (2000) and Robins, van der Vaart and Ventura (2000) and others have noted, they do not produce $p$-values that have (even asymptotically) a uniform distribution. Because output from predictive posterior model checks is not calibrated, using $p$-values based on them for model assessment is problematic.

Bayarri and Berger (2000) and Robins, van der Vaart and Ventura (2000) propose alternative distributions under which $p$-values, and thus model diagnostics, can be calculated. These include partial posterior-predictive $p$-values and conditional predictive $p$-values [Bayarri and Berger (2000)], and modifications to posterior predictive and "plug-in" $p$-values [Robins, van der Vaart and Ventura (2000)]. The attractive feature of each of these variations on more standard definitions of $p$-values is that these statistics are distributed either as $U(0,1)$ random variables, or approach $U(0,1)$ random variables as sample sizes become large. Their drawback is that they can seldom be defined and calculated in realistically complex models.

The goal of this article is to present a goodness-of-fit diagnostic that bridges the gap between diagnostics that are easy to compute but whose null distributions are unknown, and diagnostics whose null distributions are known but that cannot generally be computed. The proposed diagnostic is closely related to the classical $\chi^2$ goodness-of-fit statistic, whose properties are now briefly reviewed.

In the case of a point null hypothesis, the standard $\chi^2$ statistic may be defined as

$$R^0 = \sum_{k=1}^{K} \frac{(m_k - np_k)^2}{np_k},$$

where $m_k$ represents the number of observations observed within the $k$th partitioning element, $p_k$ the probability assigned by the null model to this interval, $K$ the number of partitions or intervals specified over the sample space and $n$ the sample size. For independent and identically distributed data satisfying certain regularity requirements, Pearson (1900) demonstrated that the asymptotic distribution of $R^0$ was $\chi^2$ with $K-1$ degrees of freedom.

The situation for composite hypotheses is more complicated. Assuming that bins are determined a priori, Cramér [(1946), pages 426–434] demonstrated that the distribution of

$$R^g = \sum_{k=1}^{K} \frac{(m_k - np_k^g)^2}{np_k^g}$$

is that of a $\chi^2$ random variable with $K - s - 1$ degrees of freedom, where $s$ denotes the dimension of the underlying parameter vector $\boldsymbol{\theta}$ and $\{p_k^g\}$ denote estimates of the bin probabilities based on either maximum likelihood estimation for the grouped data or on the minimum $\chi^2$ method. Maximum likelihood estimation for the grouped data implies maximization of the function

$$\prod_k p_k(\boldsymbol{\theta})^{m_k}$$

with respect to $\boldsymbol{\theta}$, while minimum $\chi^2$ estimation involves the determination of a value of $\boldsymbol{\theta}$ that minimizes a function related to $R^g$.

The statistic $R^g$ is the form of the $\chi^2$ test most often used in statistics, where it is routinely used to test independence in contingency tables [see, e.g., Fienberg (1980)]. In that context, grouped maximum likelihood estimation is natural. Although the Bayesian $\chi^2$ statistic proposed below can be extended for testing independence in contingency tables, this is not its intended purpose. Instead, it is intended primarily for use as a goodness-of-fit test. In this regard, the aspect of model fit assessed is similar to that examined using the classical $\chi^2$ goodness-of-fit test; namely, the proportion of counts observed in predefined parcels of the sample space is compared to the proportion of counts that are expected in these parcels under a specified probability model.

Chernoff and Lehmann (1954) considered the distribution of the $\chi^2$ statistic in the more typical situation in which values of the bin probabilities are based on maximum likelihood estimates obtained using the raw (ungrouped) data. Denote these values by $\hat{p}_k$. In this case, the distribution of the goodness-of-fit statistic is generally not one of a $\chi^2$ distribution, but instead produces a value $\hat{R}$ that has a distribution that falls stochastically between $R^0$ and $R^g$. For models containing many parameters, the gap between the degrees of freedom associated with these two statistics is large, and, as a result, the $\chi^2$ goodness-of-fit test based on the maximum likelihood estimate is usually not useful for assessing model fit in high-dimensional settings.

The goodness-of-fit statistic proposed here represents a modification of the $\chi^2$ statistics considered above. The modification, denoted by $R^B(\tilde{\boldsymbol{\theta}})$ (or more simply, by $R^B$ when no confusion arises), is obtained by fixing the values of $p_k$ and instead considering the bin counts $m_k$ as random quantities. Allocation of observations to bins is made according to the value of each observation's conditional distribution function, conditionally on a single parameter value $\tilde{\boldsymbol{\theta}}$ sampled either from the posterior distribution or the asymptotic distribution of the maximum likelihood estimator. [The statistic obtained in this way has some resemblance to the $\chi^2$ statistics considered by, e.g., Moore and Spruill (1975), although emphasis there focuses on randomized cells rather than on posterior sampling of parameter vectors.] The distinguishing feature of $R^B(\tilde{\boldsymbol{\theta}})$ is that, for many statistical models, its asymptotic distribution is $\chi^2$ on $K-1$ degrees of freedom, independently of the dimension of the parameter vector $\boldsymbol{\theta}$.

Because it is the sampling distribution of $R^B$ that has a $\chi^2$ distribution, one might argue that this procedure does not really represent a Bayesian goodness-of-fit diagnostic. However, sampling parameter values from a distribution for the purpose of inference occurs more naturally within the Bayesian paradigm, and for this reason it is likely that the proposed diagnostic will find more application there. In addition, the formal test statistics proposed below are based on the posterior distribution of $R^B$. For this reason, values of $\tilde{\boldsymbol{\theta}}$ used in the definition of $R^B$ are assumed to represent samples from the posterior distribution on the parameter vector, rather than samples generated from the asymptotic normal distribution of the maximum likelihood estimator. However, either interpretation is valid.

The remainder of the paper is organized as follows. In the next section, the Bayesian $\chi^2$ statistic $R^B$ is defined and its asymptotic properties are described. Corollaries extending these properties from i.i.d. observations to conditionally independent observations and to fixed-bin applications are presented, and strategies for combining information contained in dependent samples of $R^B$ values generated from the same posterior distribution are described. Following this, several examples that illustrate the application of this statistic and summaries from its posterior are presented. Discussion and concluding remarks appear in Section 4. Proofs of the theorem and corollaries of Section 2 appear in the Appendix.

**2. A Bayesian $\chi^2$ statistic.** To begin, let $y_1, \ldots, y_n \, (= \mathbf{y})$ denote scalar-valued, continuous, identically distributed, conditionally independent observations drawn from probability density function $f(y|\boldsymbol{\theta})$ defined with respect to Lebesgue measure and indexed by an $s$-dimensional parameter vector $\boldsymbol{\theta} \in \boldsymbol{\Theta} \subset \mathbf{R}^s$. Denote by $F(\cdot|\boldsymbol{\theta})$ and $F^{-1}(\cdot|\boldsymbol{\theta})$ the (nondegenerate) cumulative distribution and inverse distribution functions corresponding to $f(\cdot|\boldsymbol{\theta})$. To construct a sampled value $\tilde{\boldsymbol{\theta}}$ from the posterior, augment the observed sample $\mathbf{y}$ with an i.i.d. sample $v_1, \ldots, v_s$ from a $U(0,1)$ distribution. Let $p(\boldsymbol{\theta}|\mathbf{y})$ denote the posterior density of $\boldsymbol{\theta}$ based on $\mathbf{y}$, and let $p(\boldsymbol{\theta}_j|\boldsymbol{\theta}_1, \ldots, \boldsymbol{\theta}_{j-1}, \mathbf{y})$ denote the marginal conditional posterior density of $\boldsymbol{\theta}_j$ given $(\boldsymbol{\theta}_1, \ldots, \boldsymbol{\theta}_{j-1}, \mathbf{y})$. Define $\tilde{\boldsymbol{\theta}}$ implicitly by

$$
(1) \qquad v_1 = \int_{-\infty}^{\tilde{\theta}_1} p(\theta_1|\mathbf{y}) \, d\theta t_1, \ldots, v_s = \int_{-\infty}^{\tilde{\theta}_s} p(\theta_s|\tilde{\theta}_1, \ldots, \tilde{\theta}_{s-1}, \mathbf{y}) \, d\theta_s.
$$

Thus, $\tilde{\boldsymbol{\theta}}$ denotes a value of $\boldsymbol{\theta}$ sampled from the posterior distribution based on $\mathbf{y}$. Let $\boldsymbol{\theta}_0$ denote the true but unknown value of $\boldsymbol{\theta}$. The maximum likelihood estimate of $\boldsymbol{\theta}$ is denoted by $\hat{\boldsymbol{\theta}}$.

To construct the Bayesian goodness-of-fit statistic proposed here, choose quantiles $0 \equiv a_0 < a_1 < \cdots < a_{K-1} < a_K \equiv 1$, with $p_k = a_k - a_{k-1}$, $k = 1, \ldots, K$. Define $\mathbf{z}_j(\tilde{\boldsymbol{\theta}})$ to be a vector of length $K$ whose $k$th element is 0 unless

$$
(2) \qquad F(y_j|\tilde{\boldsymbol{\theta}}) \in (a_{k-1}, a_k],
$$

in which case it is 1. Finally, define

$$
\mathbf{m}(\tilde{\boldsymbol{\theta}}) = \sum_{j=1}^{n} \mathbf{z}_j(\tilde{\boldsymbol{\theta}}).
$$

It follows that the $k$th component of $\mathbf{m}(\tilde{\boldsymbol{\theta}})$, $m_k(\tilde{\boldsymbol{\theta}})$, represents the number of observations that fell into the $k$th bin, where bins are determined by the quantiles of the inverse distribution function evaluated at $\tilde{\boldsymbol{\theta}}$. Finally, define

$$
(3) \qquad R^B(\tilde{\boldsymbol{\theta}}) = \sum_{k=1}^{K} \left[ \frac{(m_k(\tilde{\boldsymbol{\theta}}) - np_k)}{\sqrt{np_k}} \right]^2.
$$

The asymptotic distribution of $R^B$ is provided in the following theorem.

THEOREM 1. *Assuming that the regularity conditions specified in the Appendix apply, $R^B$ converges to a $\chi^2$ distribution with $K - 1$ degrees of freedom as $n \to \infty$.*

The simplicity of Theorem 1 is somewhat remarkable given the complexity of the corresponding distribution of $\hat{R}$. As mentioned above, the asymptotic distribution of $\hat{R}$ does not, in general, follow a $\chi^2$ distribution. Instead, it has the distribution of the sum of a $\chi^2$ random variable with $K - s - 1$ degrees of freedom and the weighted sum of $s$ additional squared normal deviates with weights ranging from 0 to 1. In contrast, the asymptotic distribution of $R^B$ follows a $\chi^2_{K-1}$ distribution, independently of the dimension of the parameter vector $\boldsymbol{\theta}$.

Heuristically, the idea underlying Theorem 1 is that the degrees of freedom lost by substituting the grouped MLE for $\boldsymbol{\theta}$ in Pearson's $\chi^2$ statistic are exactly recovered by replacing the MLE with a sampled value from the posterior in $R^B$. That is, the $s$ degrees of freedom lost by maximizing over the grouped likelihood function to obtain $R^g$ are exactly recovered by sampling from the $s$-dimensional posterior on $\boldsymbol{\theta}$.

As a corollary, Theorem 1 can be extended to the more general case in which the functional form of the density $f(y|\boldsymbol{\theta})$ varies from observation to observation. Specifically, if the density of the $j$th observation is denoted by $f_j(y|\boldsymbol{\theta})$, with distribution and inverse distribution functions $F_j$ and $F_j^{-1}$, respectively, then the following corollary also applies.

COROLLARY 1. *Assume the conditions referenced in Theorem 1 are extended so as to provide also for the asymptotic normality of both the posterior distribution on $\boldsymbol{\theta}$ and the maximum likelihood estimator when the likelihood function is proportional to*

$$\prod_{j=1}^{n} f_j(y_j|\boldsymbol{\theta}).$$

*Assume also that the functions $f_j(\cdot|\boldsymbol{\theta})$ satisfy the same conditions implied in Theorem 1 for $f(\cdot|\boldsymbol{\theta})$. Define the $k$th component of $\mathbf{z}_j(\boldsymbol{\theta})$ to be 1 or 0 depending on whether or not*

(4) $$F_j(y_j|\tilde{\boldsymbol{\theta}}) \in (a_{k-1}, a_k],$$

*with $\mathbf{a}$ fixed. Then the asymptotic distribution of $R^B$ based on this revised definition of $\mathbf{z}_j(\boldsymbol{\theta})$ is $\chi^2$ with $K-1$ degrees of freedom.*

Outlines of the proof of Theorem 1 and the corollary appear in the Appendix.

From a practical perspective, the corollary is important because it extends the definition of $R^B$ to essentially all models in which observations are continuous and conditionally independent given the value of a finite-dimensional parameter vector.

The results cited above for continuous-valued random variables can be extended to discrete random variables in one of two ways. The most direct extension is to simply proceed as in the continuous case, using a randomization procedure to allocate counts to bins when the mass assigned to an observation spans the boundaries defining the bins. The second is to define fixed bins in the standard way based on the possible outcomes of the random variable, and to then evaluate the bin probabilities at sampled values of $\boldsymbol{\theta}$ from the posterior distribution. That is, if $f(y|\boldsymbol{\theta})$ denotes the probability mass function of a discrete random variable $y$ and

$$p_k(\tilde{\boldsymbol{\theta}}) = \frac{1}{n} \sum_{j=1}^{n} \sum_{y \in \text{bin } k} f_j(y|\tilde{\boldsymbol{\theta}}), \tag{5}$$

then the $\chi^2$ statistic $R^B$ may be redefined as

$$R^B(\tilde{\boldsymbol{\theta}}) = \sum_{k=1}^{K} \left[ \frac{(m_k - np_k(\tilde{\boldsymbol{\theta}}))}{\sqrt{np_k(\tilde{\boldsymbol{\theta}})}} \right]^2. \tag{6}$$

In this case, the asymptotic distribution of $R^B(\tilde{\boldsymbol{\theta}})$ is similar to that described above in the continuous case and is detailed in the following corollary.

COROLLARY 2. *If the regularity conditions specified in Theorem 1 apply to the discrete probability mass function $f(y|\boldsymbol{\theta})$, then, using predefined bins and the definition of the bin probabilities given in (5), the distribution of $R^B(\tilde{\boldsymbol{\theta}})$ as defined in (6) converges to a $\chi^2$ distribution with $K-1$ degrees of freedom as $n \to \infty$.*

The asymptotic $\chi^2$ distribution of $R^B(\tilde{\boldsymbol{\theta}})$ described in the theorem and corollaries above is achieved when a large sample of independent observations is drawn from a sampling density, and a value of $\tilde{\boldsymbol{\theta}}$ is drawn from the posterior induced by this observation. However, when two values of $\tilde{\boldsymbol{\theta}}$ are drawn from the same posterior distribution (i.e., based on the same observation), the values of $R^B$ that result are correlated. This correlation complicates the interpretation of test statistics defined with respect to posterior distribution on $R^B$ values.

Combining information across a posterior sample of $R^B$ values might be accomplished in a variety of ways, including modifications of the methodologies proposed in Verdinelli and Wasserman (1998) or Robert and Rousseau (2002). Another possibility is to simply report the proportion of $R^B$ values drawn from the posterior distribution that exceeds a specified critical value from their nominal $\chi^2_{K-1}$ distribution. For a given data vector and probability model, such a procedure might lead to a statement that, say, 90% of $R^B$ values generated from the posterior distribution exceeded the 95th percentile of the reference $\chi^2$ distribution.

Though decidedly non-Bayesian, such a report is convenient from several perspectives. By reporting the proportion of $R^B$ values that exceeds the critical value of the test, the unpalatable aspect of basing a goodness-of-fit test on a randomly selected value of $R^B$ is avoided. It is also straightforward to compare the proportion of $R^B$ values that exceeds the critical value of the test to the size of the test; if the $R^B$ values did represent independent draws from their nominal $\chi^2$ distribution, the proportion of values falling in the critical region of the test would exactly equal the size of the test. Any excess in this proportion must therefore be attributed either to dependence between the sampled values of $R^B$ from the given posterior or lack of fit. Finally, and perhaps most importantly, this strategy requires almost no computational effort. In most practical Bayesian models, values of $R^B$ can be computed almost as an afterthought within the MCMC schemes used to sample from the posterior distribution of the parameter vector.

In the event that formal significance tests must be performed to assess model adequacy, they can be based on a comparison of the observed value of a summary statistic based on the posterior distribution of $R^B$ values to an approximation of the sampling distribution of the summary statistic induced by repeated sampling of the data vector. The summary statistic considered here is defined as the posterior probability that a value of $R^B$ drawn from the posterior distribution (based on a single value of **y**) exceeds the value of a $\chi^2_{K-1}$ random variable. This probability, denoted by $A$, is related to a commonly used quantity in signal detection theory and represents the area under the ROC curve [e.g., Hanley and McNeil (1982)] for comparing the joint posterior distribution of $R^B$ values to a $\chi^2_{K-1}$ random variable. The expected value of $A$, if taken with respect to the joint sampling distribution of **y** and the posterior distribution of $\boldsymbol{\theta}$ given **y**, would be 0.5. Large deviations in the expected value of $A$ from 0.5, when the expectation is taken with respect to the posterior distribution of $\boldsymbol{\theta}$ for a fixed value of **y**, indicate model lack of fit.

Unfortunately, approximating the sampling distribution of $A$ is a numerically burdensome endeavor, and calculating it obviates many of the advantages that are gained by using a test statistic with a known reference distribution. To a large extent, the computations required to approximate $A$'s sampling distribution are as complicated as, or even more complicated than, similar techniques used to approximate the sampling distribution of discrepancy functions used in posterior-predictive model checks [e.g., Sinharay and Stern (2003)]. However, knowing the nominal value of $A$ makes this computation unnecessary when the observed value of $A$ falls within several hundredths of 0.5 or is smaller than 0.5. Procedures for approximating the sampling distribution of $A$ for the purpose of determining the significance of departures of the observed value of $A$ from 0.5 are described in the examples using methodology delineated by Dey, Gelfand, Swartz and Vlachos (1998).

As an aside, it is interesting to compare the test statistic $R^B$ and its reference distribution to the $\chi^2$ discrepancy function and its reference distribution as proposed in Gelman, Meng and Stern (1996). The reference distribution of $R^B(\tilde{\boldsymbol{\theta}})$ is obtained by sampling **y** from its "true" distribution $F(\cdot|\boldsymbol{\theta}_0)$, and then sampling a single value of $\tilde{\boldsymbol{\theta}}$ from the posterior distribution of $\boldsymbol{\theta}$ given $\mathbf{y}$. The resulting distribution is asymptotically $\chi^2_{K-1}$; this result is unrelated to posterior-predictive distributions or samples drawn from them. In contrast, the reference distribution of the $\chi^2$ discrepancy function proposed by Gelman, Meng and Stern (1996) is obtained as the distribution of the statistic

$$(7) \qquad \sum_{i=1}^n \frac{(y_i^{pp} - E(y_i^{pp}|\boldsymbol{\theta}))^2}{\operatorname{Var}(y_i^{pp}|\boldsymbol{\theta})}$$

induced by repeatedly drawing values $\mathbf{y}^{pp} = (y_1^{pp}, \ldots, y_n^{pp})$ from the posterior-predictive density based on the observed data vector $\mathbf{y}$. As Gelman, Meng and Stern point out, this statistic does not have a $\chi^2$ distribution.

The power characteristics of the Bayesian $\chi^2$ statistics defined above, like their classical counterparts, depend on the selection of the bin probabilities $p_k$. Clearly, consistency of derived tests against general alternatives requires that $K \to \infty$ as $n \to \infty$. On the other hand, as many authors have noted [see, e.g., Koehler and Gan (1990) for a review of this topic], using too many cells can result in a significant loss of power.

A general criterion for choosing cell probabilities was proposed by Mann and Wald (1942), who suggested the use of $3.8(n-1)^{0.4}$ equiprobable cells. Subsequent authors [e.g., Williams (1950), Watson (1957), Hamdan (1963), Dahiya and Gurland (1973), Gvanceladze and Chibisov (1979), Best and Rayner (1981), Quine and Robinson (1985) and Koehler and Gan (1990)] found that the Mann–Wald criterion often results in too many bins and loss of power. Based on numerical simulations of seven classes of alternative probability models, Koehler and Gan (1990) noted that near-optimal power against a Gaussian null model was obtained when the Mann–Wald criterion was divided by 4. Such a rule also finds approximate agreement with simulation results reported by Kallenberg, Oosterhoff and Schriever (1985) (although they also recommend the use of nonequiprobable cells against certain types of alternative hypotheses). This rule of thumb, which may be approximately reformulated as taking $n^{0.4}$ equiprobable cells, was found to yield nearly optimal results in the examples described below.

## 3. Examples.

3.1. *Goodness-of-fit tests under a normal model with unknown mean and variance.* In this example, the distribution of $R^B$ under a normal model is investigated and compared with the distributions of $\hat{R}$ and $R^g$. Posterior samples of $R^B$ generated from a single data vector are used in ROC-type analyses to generate a summary model diagnostic. The power of this test statistic is investigated and compared to the power of the test statistic $R^g$ when data are generated under nonnormal alternatives.

Let $\mathbf{y} = (y_1, \ldots, y_{50})$ denote a random sample from a standard normal distribution. For purposes of illustration, assume that the mean $\mu$ and variance $\sigma^2$ of the data are unknown and that the joint prior assumed for $(\mu, \sigma)$ is proportional to $1/\sigma$. Let $(\tilde{\mu}, \tilde{\sigma})$ denote a sampled value from the posterior distribution based on $\mathbf{y}$.

For a given data vector $\mathbf{y}$ and posterior sample $(\tilde{\mu}, \tilde{\sigma})$, bin counts $m_k(\tilde{\mu}, \tilde{\sigma})$ are determined by counting the number of observations $y_i$ that fall into the interval $(\tilde{\sigma}\Phi^{-1}(a_{k-1}) + \tilde{\mu}, \tilde{\sigma}\Phi^{-1}(a_k) + \tilde{\mu})$, where $\Phi^{-1}(\cdot)$ denotes the standard normal quantile function. Based on these counts, $R^B(\tilde{\mu}, \tilde{\sigma})$ is calculated according to (3).

Figure 1 depicts a quantile-quantile plot of $R^B$ values calculated for 10,000 independent samples of $\mathbf{y}$. Each value of $R^B$ depicted in this plot corresponds to a single draw of $(\mu, \sigma)$ from the posterior distribution based on a single observation vector $\mathbf{y}$. In accordance with the rule-of-thumb discussed in Section 2, five equiprobable bins were used in the definition of $R^B$. As expected, the distribution of $R^B$ closely mimics that of a $\chi_4^2$ random variable.

The normal deviates used in the construction of Figure 1 were also used to compute the classical $\chi^2$ statistic based on the maximum likelihood estimates of $\mu$ and $\sigma$ (i.e., using the ungrouped data). The quantile-quantile plot of 10,000 $\hat{R}$ values obtained from these data is displayed in Figure 2. In this plot, the $\hat{R}$ values have been plotted against the expected order statistics from a $\chi_2^2$ random variable. Five equal probability bins based on the standard normal distribution were also used to define these $\hat{R}$ values. As might be expected, the plotted $\chi^2$ values display some deviation from their approximate $\chi_2^2$ distribution.

Grouped maximum likelihood estimates were also used to calculate $R^g$ values using these normal samples. The corresponding quantile-quantile plot for the 10,000 $R^g$ values is displayed in Figure 3; as expected, these values demonstrate substantially better agreement with a $\chi_2^2$ random variable than do the values depicted in Figure 2.

Returning to the investigation of the properties of $R^B$, Figure 1 demonstrates excellent agreement between this statistic and its asymptotic distribution. To illustrate its power in detecting departures from the normal model, suppose now that the experiment above is repeated with independent Student $t$ variates substituted for the normal deviates. That is, the actual observation vectors used in the simulation represent Student $t$ variates, but the statistical model used to calculate values of $R^B$ is still based on the assumption that the data are normally distributed. The degrees of freedom of the $t$ variates used in this experiment range from 1 to 10, and for each value within this range, 10,000 independent samples of size 50 were drawn.

To study the power of the statistic $R^B$ in detecting departures from normality in this experiment, formal significance tests were performed using the

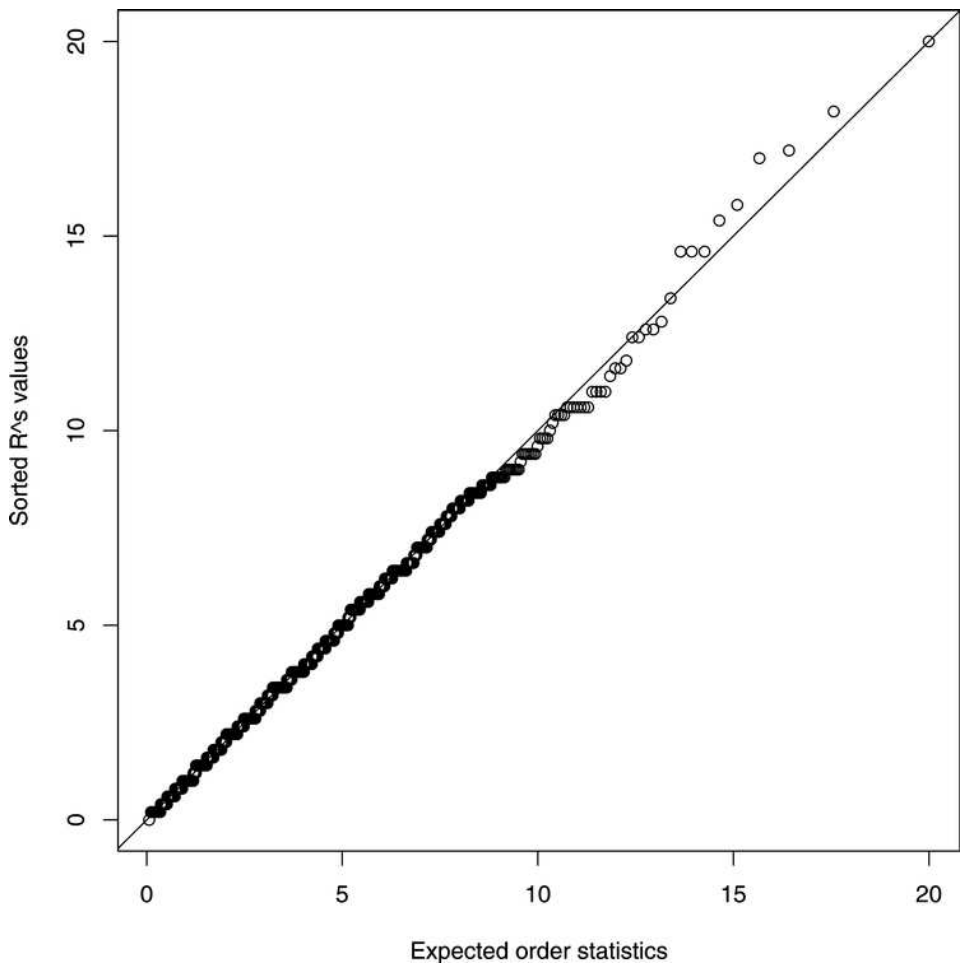

FIG. 1. *Quantile-quantile plot of $R^B$ values for i.i.d. normal data. Values of $R^B$ displayed in this plot were determined from independent samples of* 50 *standard normal deviates, and are plotted against the expected order statistics from a $\chi_4^2$ distribution. Posterior samples of the mean and variance were estimated using reference priors and observations were binned into five bins of equal probability [i.e., $\mathbf{a} = (0, 0.2, 0.4, 0.6, 0.8, 1)$].*

statistic $A$ described in Section 2. This statistic may be defined formally as

(8)  $$A = \mathbf{Pr}_{\tilde{\boldsymbol{\theta}}|\mathbf{y}}(R^B(\tilde{\boldsymbol{\theta}}) > X), \qquad X \sim \chi^2_{K-1},$$

and, in repeated sampling of both $\mathbf{y}$ and $\boldsymbol{\theta}$ given $\mathbf{y}$, has a nominal value of 0.5. Numerically, the value of $A$, for a fixed data vector $\mathbf{y}$, can be approximated in a straightforward way using Monte Carlo integration.

Formal model assessment using the statistic $A$ can be based on approximating the sampling distribution of $A$ using "posterior-predictive-posterior" model checks [e.g., Dey, Gelfand, Swartz and Vlachos (1998)]. That is, sampled values $\tilde{\boldsymbol{\theta}}$ from the posterior can be used to generate posterior-predictive observations $\mathbf{y}^{pp}$ according to $f(\cdot|\tilde{\boldsymbol{\theta}})$. In large samples, values of $\tilde{\boldsymbol{\theta}}$ will be close to $\boldsymbol{\theta}_0$, and so the distribution of $\mathbf{y}^{pp}$ will be close to the distribution of $\mathbf{y}$. Posterior-predictive-posterior values of $A^{pp}$ can be generated for each value of $\mathbf{y}^{pp}$ by averaging $R^B$, computed from $\mathbf{y}^{pp}$, over the posterior distribution of $\boldsymbol{\theta}$ induced by $\mathbf{y}^{pp}$. Values of $A^{pp}$ obtained from this procedure approximate the sampling variability of the summary test statistic $A$ that can be attributed to computing the probability in (8) using the posterior distribution of $\boldsymbol{\theta}$ for a given value of $\mathbf{y}$, without averaging over the distribution of $\mathbf{y}$. The value of $A$ obtained for the original data vector can then be compared to the empirical distribution of the values of $A^{pp}$ obtained from the posterior distribution on the posterior-predictive data.

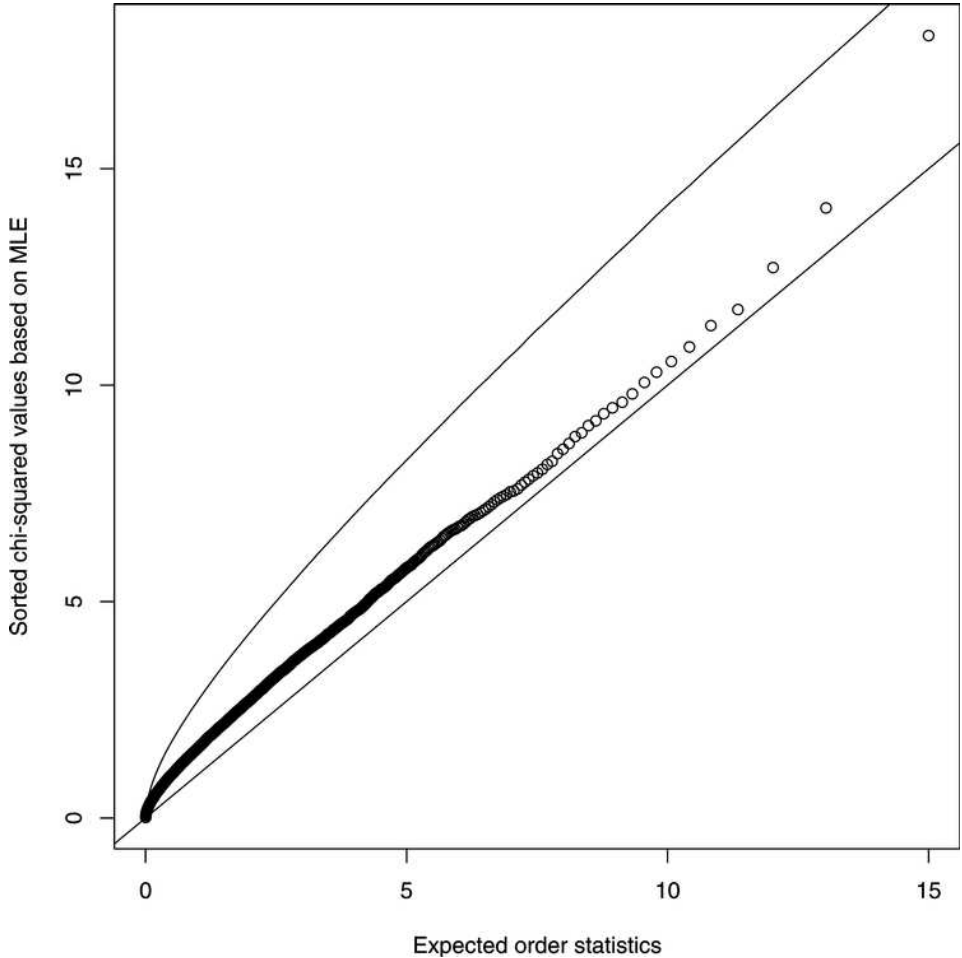

FIG. 2. *Quantile-quantile plot of $\hat{R}$ values for i.i.d. normal data. Values of $\hat{R}$ displayed in this plot were each determined from a separate sample of 50 standard normal deviates, and are plotted against the expected order statistics from a $\chi^2_2$ distribution. For comparison, the top curve depicts values of expected order statistics from a $\chi^2_4$ distribution.*

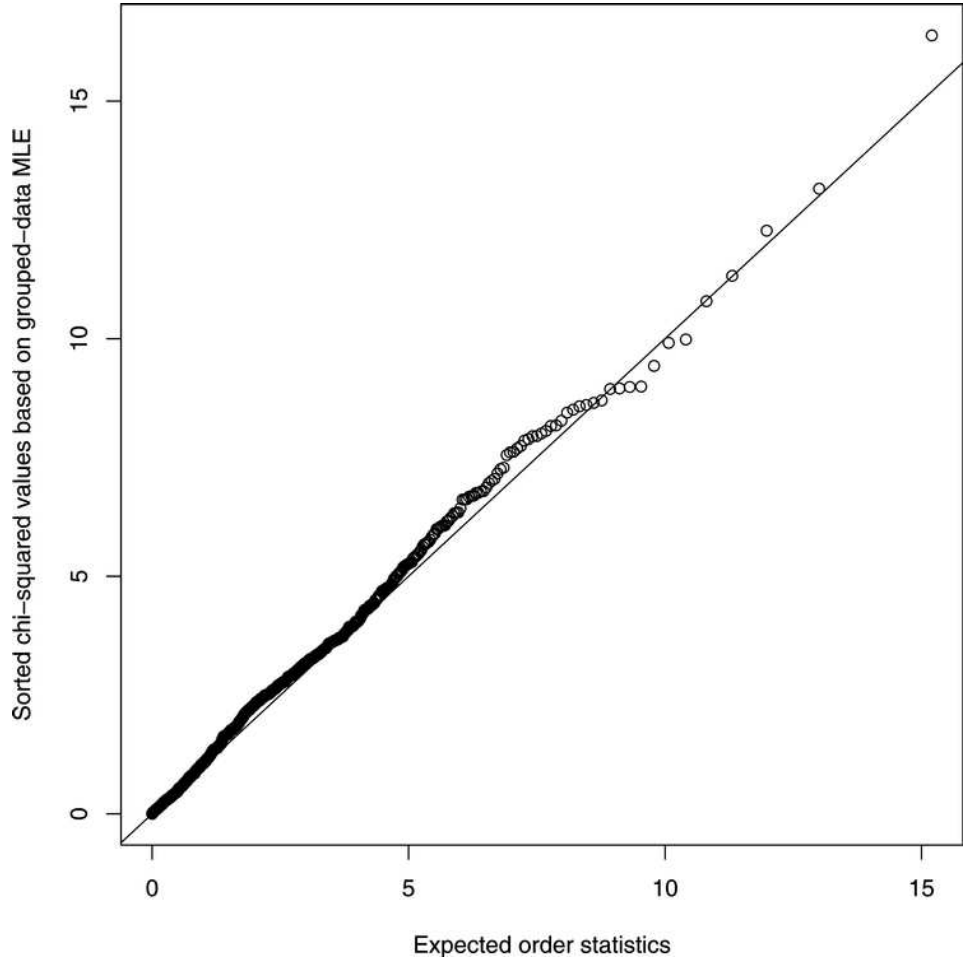

FIG. 3. *Quantile-quantile plot of $R^g$ values for i.i.d. normal data. Values of $R^g$ displayed in this plot were each determined from a separate sample of* 50 *standard normal deviates, and are plotted against the expected order statistics from their asymptotic $\chi_2^2$ distribution.*

In principle, exactly this procedure can be implemented to calculate the probability that the test statistic $A$, based on a random sample of $t$ variates, falls into the critical region of a test based on the empirical distribution of sampled values $A^{pp}$. In this case, however, it is not necessary to generate values of $A^{pp}$ for each sample of $t$ variates. Under the normal model, values of $R^B$ are invariant to shifts in location and scale of the data, so the sampling distribution of $A$, for any future draw of 50 i.i.d. normal deviates, can be approximated by the empirical distribution of $A$ values obtained under the normal sampling scheme used at the beginning of this example. It follows that critical regions for significance tests based on $A$ are exact under this model, save for the Monte Carlo error encountered in the empirical approximation of their distribution.

Figure 4 displays the proportion of times in 10,000 draws of $t$ samples that the value of the test statistic $A$ was larger than the 0.95 quantile of the sampled values of $A^{pp}$. For comparison, the observed power of the test based on the grouped-maximum-likelihood $\chi^2$ statistic $R^g$ at the 5% level is also shown, as is the observed power obtained using a randomized test based on only a single value of $R^B$. To facilitate comparison with the distribution of $R^B$, five equiprobable bins from a standard normal distribution were used in the definition of $R^g$.

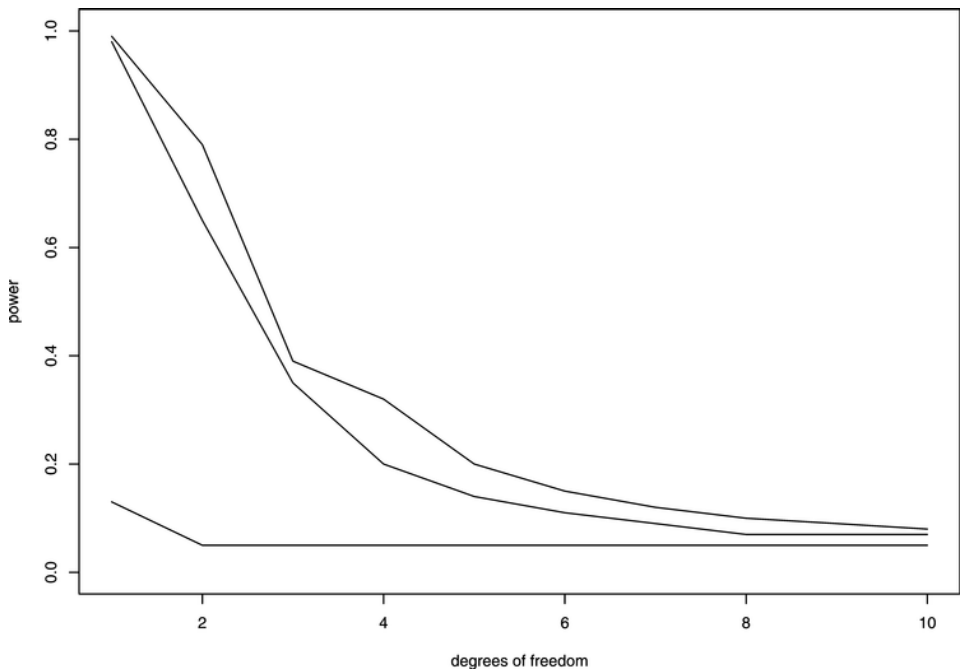

Fig. 4. *Power of test statistics $A$, $R^B$ and $R^g$ in detecting departures from normality when data are distributed according to $t$ distributions. The uppermost curve depicts the power of the test statistic $A$ against $t$ alternatives with degrees of freedom displayed on the horizontal axis. The curve in the middle depicts the corresponding power of a single value of $R^B$ when compared to a $\chi_4^2$ distribution. The curve at the bottom of the plot represents the power of $R^g$ against the $t$ alternatives. All values of the power refer to the power of the test statistics in rejecting the null hypothesis of normality in significance tests of size 0.05 and samples of size* 50.

From Figure 4, it is clear that the test statistic $A$ offers substantially better power than $R^g$ against this class of alternative models. Part of this advantage stems from the symmetry and unimodality of the alternative hypotheses, which $R^g$ is ill-equipped to accommodate, and part from the fact that the bins used in the definition of $R^g$ were fixed according to the null hypothesis. Substantially better power could be obtained by using the test statistic $\hat{R}$ with bins based on the particular **y** vector observed, but such tests do not achieve their nominal levels of significance. Perhaps surprisingly, the power of a randomized test based on a single value of $R^B$ is comparable to the power of $A$ based on the complete posterior distribution of $R^B$ values.

3.2. *Lip cancer data.* Spiegelhalter, Best, Carlin and van der Linde (2002) describe a re-analysis of lip cancer incidence data originally considered by Clayton and Kaldor (1987). Their purpose in examining these data was to illustrate the use of the deviance information criterion (DIC) to select from among five potential models for the number of lip cancer cases, $y_i$, observed in 56 Scottish districts as a function of available age and sex adjusted expected rates $E_i$. These data and models are reconsidered here for the related purpose of assessing which of the models provides an adequate probabilistic description of the data.

Following the Spiegelhalter et al. analysis of these data, begin by assuming that $y_i$ is Poisson with mean $\mu_i = \exp(\theta_i)E_i$. Five models for $\theta_i$ are considered:

1. $\theta_i = \alpha_0$, $\alpha_0$ a constant.
2. $\theta_i = \alpha_0 + \gamma_i$, $\gamma_i$ exchangeable random effects.
3. $\theta_i = \alpha_0 + \delta_i$, $\delta_i$ spatial random effects with a conditional autoregressive prior [e.g., Besag (1974)].
4. $\theta_i = \alpha_0 + \delta_i + \gamma_i$, $\delta_i$ and $\gamma_i$ as above.
5. $\theta_i = \alpha_i$, $\alpha_i$ uniform on $(-\infty, \infty)$.

Five thousand, thinned posterior samples of $\boldsymbol{\mu} = \{\mu_i\}$ were generated for each of these models using WinBUGS code [Spiegelhalter, Thomas and Best (2000)] kindly provided by Dr. Best. For each sampled value of $\mu_i$, the Poisson counts $y_i$ were assigned to one of five equiprobable bins defined according to the Poisson distribution function evaluated at $y_i$ for the given value of $\mu_i$. In those cases for which the probability mass function assigned to $y_i$ spanned more than one bin, allocation to a single bin was performed randomly according to the proportion of mass assigned to the bins. Averaging over all posterior samples of $\boldsymbol{\mu}$ for a given model yielded the values of $A$ depicted in Table 1. Because 56 data points were available, five bins were again used in the definition of the individual values of $R^B$. The proportion of $R^B$ values exceeding the 95th quantile from a $\chi_4^2$ distribution was computed using the

TABLE 1
*Values of the goodness-of-fit statistic A and the proportion of critical $R^B$ values for models of lip cancer incidence data*

| Model | $A$ | Proportion of $R^B > 9.49$ | DIC |
|---|---|---|---|
| 1 | 0.999 | 1.000 | 382.7 |
| 2 | 0.517 | 0.055 | 104.0 |
| 3 | 0.538 | 0.076 | 89.9 |
| 4 | 0.537 | 0.075 | 89.7 |
| 5 | 0.677 | 0.198 | 111.7 |

The second column provides the value of the summary statistic $A$ achieved for each model. The third column lists the proportion of posterior samples of $R^B$ that exceeded the 95th quantile of a $\chi_4^2$ distribution for each model. DIC values obtained under the "mean" parameterization are listed for comparison.

posterior sample $\boldsymbol{\mu}$. No posterior-predictive or posterior-predictive-posterior computations were performed to obtain these values.

In Table 1, both the large value of $A$ and the large proportion of $R^B$ values exceeding the 95th quantile of the $\chi_4^2$ distribution provide a clear indication of lack of fit for the first model. Lack of fit in this instance can be attributed to the failure of the model to adjust for district effects; the posterior mean of the number of counts assigned to the five bins was $(16.0, 4.9, 5.2, 7.1, 22.8)$.

The values of $A$ and proportions of extreme values of $R^B$ reported in rows 2–4 do not suggest lack of fit of the aspect of these models being tested by the $\chi^2$ test.

The most interesting row in Table 1 is the last, which corresponds to fitting a separate Poisson model for each observation. The value of $A$ for this model is 0.68, and nearly 20% of $R^B$ values generated from its posterior—nearly four times the number expected—exceeded the 5% critical value from the $\chi_4^2$ distribution.

At first glance, one might suspect that these suspicious values arise from overfitting. However, the last model generates the most dispersed posterior distribution of any of the models considered, since only one observation figures into the marginal posterior of each $\mu_i$. Instead, the difficulty with this model arises from the prior assumptions made on $\boldsymbol{\mu}$. The assumption of a uniform prior on $\theta_i$ implies a prior for the mean of each Poisson observation proportional to $1/\mu_i$; this prior shrinks the estimate of every $\mu_i$ toward 0. This results in an overabundance of counts in the higher bins and larger than expected values of $R^B$. The posterior mean of the bin counts for this model was $(8.4, 9.8, 10.9, 12.1, 14.8)$. Refitting the fifth model with noninformative priors proportional to $1/\sqrt{\mu_i}$ yielded a value of $A = 0.501$ and only 4.7% of $R^B$ values exceeding 9.49.

It is also interesting to compare the values in the second and third columns of this table with those provided for the DIC. All statistics suggest inadequacy of the first model, though for different reasons. For the first model, the high values of $A$ and $R^B$ suggest that the data do not follow Poisson distributions with a common scaling of adjusted expected rates. The value of the DIC statistic suggests either that the model does not fit the data or that it is not as precise in predicting the data as the other models considered. An advantage of the $\chi^2$ statistics in this case is that their values are interpretable without fitting alternative models.

The comparatively large value of the DIC statistic for the second model can be attributed to greater dispersion in its posterior as compared to posterior dispersion of the third and fourth models, even though the exchangeable model appears to adequately represent variation in the observed data. The comparatively large value of DIC reported for the fifth model reflects some combination of lack of fit and a posterior that is more dispersed than others considered.

**4. Extensions.** In addition to providing a convenient mechanism for assessing model adequacy, values of $R^B$ generated from a posterior distribution may prove useful both as a convergence diagnostic for MCMC algorithms and for detecting errors written in computer code to implement these algorithms.

Monitoring values of $R^B$ generated within an MCMC algorithm provides a rudimentary convergence diagnostic for slow-mixing chains. In fact, exceedances of $R^B$ over a prespecified quantile from its null distribution can be incorporated formally into the convergence diagnostics proposed in Raftery and Lewis ([1992](#)). To the extent that such exceedances are adequately described by a two-state Markov chain, the use of $R^B$ in this context eliminates the requirement to assess convergence on a parameter-by-parameter basis, as is normally done in Raftery and Lewis' diagnostic scheme. It also provides a natural mechanism for determining whether burn-in has occurred.

A less obvious but perhaps equally important use of the $R^B$ statistic involves code verification. Many practitioners currently fit models using customized code written for their specific application, a practice that frequently results in coding errors that are difficult to detect. This problem can be largely overcome by simply monitoring the distribution of $R^B$, which, in my experience, tends to deviate substantially from its null distribution when a model has been misspecified or miscoded.

**5. Discussion.** Goodness-of-fit tests based on the statistic $R^B$ provide a simple way of assessing the adequacy of model fit in many Bayesian models. Essentially, the only requirement for their use is that observations be conditionally independent. From a computational perspective, such statistics can be calculated in a straightforward way using output from existing MCMC algorithms.

Approximating the sampling distribution of $A$, though conceptually straightforward, does introduce an additional computational burden, but is necessary only when the achieved value of $A$ is "significantly" larger than 0.5. Significance of $A$ in this context has a natural interpretation in terms of the posterior probability that a sampled value of $R^B$ exceeds a random variable drawn from its nominal $\chi^2$ distribution. In this regard, values of $A$ that are close to 0.5 may indicate adequate model fit for the purposes of a given analysis even when the sampling distribution of $A^{pp}$ would permit rejection of the model in a significance test.

Aside from applications in Bayesian model assessment, the $\chi^2$ statistic proposed here can be extended, albeit somewhat awkwardly, to models estimated using maximum likelihood. In that setting, parameter values can be sampled from their asymptotic normal distribution and used as if they were sampled from a posterior distribution. Although not entirely palatable from a classical perspective, such a procedure does provide a mechanism for conducting a (suboptimal) goodness-of-fit test for complicated models in which alternative tests may be difficult to perform.

## APPENDIX

**Outlines of proofs of theorems and corollaries.** The proofs of Theorem 1 and Corollary 1 are based largely on the proof given in Chernoff and Lehmann (1954) in establishing the asymptotic distribution of $\hat{R}$.

Assume that conditions specified in Cramér [(1946), pages 426 and 427] and Chen (1985) apply. Cramér specifies conditions that are sufficient for establishing the distribution of the $\chi^2$ goodness-of-fit statistic when evaluated at the parameter vector maximizing the likelihood estimate based on the grouped data, whereas Chen's conditions are sufficient for establishing the asymptotic normality of the posterior distribution. Essentially, these conditions require that the likelihood be a smooth function of the parameter vector $\boldsymbol{\theta}$ in an open interval containing $\boldsymbol{\theta}_0$ (the true value of $\boldsymbol{\theta}$), that the posterior distribution concentrate around a point in this interval, that the information contained in the observations increase with sample size and that the prior assign nonnegligible mass to the interval containing $\boldsymbol{\theta}$. In addition, assume that all third-order partial derivatives of $f(y|\boldsymbol{\theta})$ [or, in the case of the corollary, $f_j(y|\boldsymbol{\theta})$] with respect to the components of $\boldsymbol{\theta}$ exist and are bounded in an open interval containing $\boldsymbol{\theta}_0$. Finally, note that all expectations and statements regarding probabilistic orders of magnitude described below are computed under the sampling distribution of $\mathbf{y}$ given $\boldsymbol{\theta}_0$.

The following lemmas are needed.

LEMMA A.1. *Under the conditions stated above, if $\hat{\boldsymbol{\theta}}$ refers to the maximum likelihood estimate of $\boldsymbol{\theta}$, $\tilde{\boldsymbol{\theta}}$ refers to a value of $\boldsymbol{\theta}$ sampled from the posterior distribution and $m_k(\cdot)$ refers to the number of counts assigned to the $k$th bin at a specified value of $\boldsymbol{\theta}$, then*

$$\text{(9)} \qquad \frac{1}{\sqrt{n}}(m_k(\tilde{\boldsymbol{\theta}}) - m_k(\hat{\boldsymbol{\theta}})) = \frac{1}{\sqrt{n}}(m_k^*(\tilde{\boldsymbol{\theta}}) - m_k^*(\hat{\boldsymbol{\theta}})) + o_p(1)$$

$$\text{(10)} \qquad = \frac{1}{\sqrt{n}} \sum_{i=1}^{s} \frac{\partial m_k^*(\hat{\boldsymbol{\theta}})}{\partial \theta_i}(\tilde{\theta}_i - \hat{\theta}_i) + o_p(1),$$

*where*

$$m_k^*(\boldsymbol{\theta}) = n\mathbf{E}[\text{Ind}(y \in [F^{-1}(a_{k-1}|\boldsymbol{\theta}), F^{-1}(a_k|\boldsymbol{\theta})])].$$

PROOF. Expanding $m_k^*(\tilde{\boldsymbol{\theta}})$ in a Taylor series expansion about $m_k^*(\hat{\boldsymbol{\theta}})$ yields

$$\text{(11)} \qquad m_k^*(\tilde{\boldsymbol{\theta}}) - m_k^*(\hat{\boldsymbol{\theta}}) = \sum_{i=1}^{s} \frac{\partial m_k^*(\hat{\boldsymbol{\theta}})}{\partial \theta_i}(\tilde{\theta}_i - \hat{\theta}_i) + O_p\left(\frac{1}{n}\right).$$

Define

$$\Delta z_{k,j} = z_{k,j}(\tilde{\boldsymbol{\theta}}) - z_{k,j}(\hat{\boldsymbol{\theta}}).$$

Then

$$|\Delta z_{k,j}| \leq \text{Ind}(y_j \in [\min(F^{-1}(a_{k-1}|\tilde{\boldsymbol{\theta}}), F^{-1}(a_{k-1}|\hat{\boldsymbol{\theta}})),$$
$$\max(F^{-1}(a_{k-1}|\tilde{\boldsymbol{\theta}}), F^{-1}(a_{k-1}|\hat{\boldsymbol{\theta}}))])$$
$$+ \text{Ind}(y_j \in [\min(F^{-1}(a_k|\tilde{\boldsymbol{\theta}}), F^{-1}(a_k|\hat{\boldsymbol{\theta}})),$$
$$\max(F^{-1}(a_k|\tilde{\boldsymbol{\theta}}), F^{-1}(a_k|\hat{\boldsymbol{\theta}}))]).$$

Because $(\hat{\boldsymbol{\theta}} - \tilde{\boldsymbol{\theta}})$ is $O_p(1/\sqrt{n})$, $\sqrt{n}\Delta z_{k,j} = O_p(1)$. It follows that

$$\sqrt{n} \sum_j \Delta z_{k,j}/n = \frac{m_k(\tilde{\boldsymbol{\theta}}) - m_k(\hat{\boldsymbol{\theta}})}{\sqrt{n}} = \frac{m_k^*(\tilde{\boldsymbol{\theta}}) - m_k^*(\hat{\boldsymbol{\theta}})}{\sqrt{n}} + o_p(1).$$

Substituting this expression into (11) yields (10). □

COROLLARY A.2. *The previous lemma also applies if $\boldsymbol{\theta}_0$ is substituted for $\tilde{\boldsymbol{\theta}}$, that is,*

$$\frac{1}{\sqrt{n}}(m_k(\boldsymbol{\theta}_0) - m_k(\hat{\boldsymbol{\theta}})) = \frac{1}{\sqrt{n}}(m_k^*(\boldsymbol{\theta}_0) - m_k^*(\hat{\boldsymbol{\theta}})) + o_p(1)$$

$$= \frac{1}{\sqrt{n}} \sum_{i=1}^{s} \frac{\partial m_k^*(\hat{\boldsymbol{\theta}})}{\partial \theta_i}(\theta_{0,i} - \hat{\theta}_i) + o_p(1).$$

LEMMA A.3. *Define*

$$\hat{p}_k = F[F^{-1}(a_k|\boldsymbol{\theta}_0)|\hat{\boldsymbol{\theta}}] - F[F^{-1}(a_{k-1}|\boldsymbol{\theta}_0)|\hat{\boldsymbol{\theta}}] = \int_{F^{-1}(a_{k-1}|\boldsymbol{\theta}_0)}^{F^{-1}(a_k|\boldsymbol{\theta}_0)} f(y|\hat{\boldsymbol{\theta}})\,dy. \tag{12}$$

*Then, under the conditions stated above,*

$$\hat{p}_k - p_k = \frac{1}{n}(m_k^*(\boldsymbol{\theta}_0) - m_k^*(\hat{\boldsymbol{\theta}})) + O_p\!\left(\frac{1}{n}\right). \tag{13}$$

PROOF. For notational simplicity, define

$$G(\gamma, \delta; c) = F[F^{-1}(c|\gamma)|\delta]$$

and

$$H_i(\gamma; c) = \left.\frac{\partial G(\gamma, \delta; c)}{\partial \delta_i}\right|_{\delta=\gamma}.$$

Then, noting that $m_k^*(\boldsymbol{\theta}_0) = np_k = n(G(\boldsymbol{\theta}_0, \boldsymbol{\theta}_0, a_k) - G(\boldsymbol{\theta}_0, \boldsymbol{\theta}_0, a_{k-1}))$,

$$(\hat{p}_k - p_k) - \frac{1}{n}(m_k^*(\boldsymbol{\theta}_0) - m_k^*(\hat{\boldsymbol{\theta}}))$$

$$= [G(\boldsymbol{\theta}_0, \hat{\boldsymbol{\theta}}; a_k) - G(\boldsymbol{\theta}_0, \hat{\boldsymbol{\theta}}; a_{k-1})] + [G(\hat{\boldsymbol{\theta}}, \boldsymbol{\theta}_0; a_k) - G(\hat{\boldsymbol{\theta}}, \boldsymbol{\theta}_0; a_{k-1})] - 2p_k$$

$$= \left[\sum_i H_i(\boldsymbol{\theta}_0; a_k)(\hat{\theta}_i - \boldsymbol{\theta}_{0,i}) - \sum_i H_i(\boldsymbol{\theta}_0; a_{k-1})(\hat{\theta}_i - \boldsymbol{\theta}_{0,i})\right]$$

$$+ \left[\sum_i H_i(\hat{\boldsymbol{\theta}}; a_k)(\boldsymbol{\theta}_{0,i} - \hat{\theta}_i) - \sum_i H_i(\hat{\boldsymbol{\theta}}; a_{k-1})(\boldsymbol{\theta}_{0,i} - \hat{\theta}_i)\right] + O_p\!\left(\frac{1}{n}\right)$$

$$= \sum_i [H_i(\boldsymbol{\theta}_0; a_k) - H_i(\hat{\boldsymbol{\theta}}; a_k)](\hat{\theta}_i - \boldsymbol{\theta}_{0,i})$$

$$- \sum_i [H_i(\boldsymbol{\theta}_0; a_{k-1}) - H_i(\hat{\boldsymbol{\theta}}; a_{k-1})](\hat{\theta}_i - \boldsymbol{\theta}_{0,i}) + O_p\!\left(\frac{1}{n}\right)$$

$$= \sum_h \sum_i \left[\frac{\partial H_i(\boldsymbol{\theta}_0; a_k)}{\partial \boldsymbol{\theta}_{0,h}} - \frac{\partial H_i(\boldsymbol{\theta}_0; a_{k-1})}{\partial \boldsymbol{\theta}_{0,h}}\right](\hat{\theta}_h - \boldsymbol{\theta}_{0,h})(\hat{\theta}_i - \boldsymbol{\theta}_{0,i}) + O_p\!\left(\frac{1}{n}\right)$$

$$= O_p\!\left(\frac{1}{n}\right). \qquad \square$$

COROLLARY A.4.

$$\sqrt{n}(\hat{p}_k - p_k) = \frac{1}{\sqrt{n}}(m_k(\boldsymbol{\theta}_0) - m_k(\hat{\boldsymbol{\theta}})) + O_p\!\left(\frac{1}{\sqrt{n}}\right).$$

PROOF OF THEOREM 1. Decompose the terms appearing in (3) as follows:

$$(14) \quad \frac{m_k(\tilde{\boldsymbol{\theta}}) - np_k}{\sqrt{np_k}} = \frac{m_k(\tilde{\boldsymbol{\theta}}) - m_k(\hat{\boldsymbol{\theta}})}{\sqrt{np_k}} - \frac{m_k(\boldsymbol{\theta}_0) - m_k(\hat{\boldsymbol{\theta}})}{\sqrt{np_k}} + \frac{m_k(\boldsymbol{\theta}_0) - np_k}{\sqrt{np_k}}.$$

From Lemma A.1 and Corollary A.2, the first two terms on the right-hand side of (14) are asymptotically equivalent to

$$(15) \quad \frac{\sum_i \partial m_k^*(\hat{\boldsymbol{\theta}})/\partial \theta_i (\tilde{\theta}_i - \hat{\theta}_i)}{\sqrt{np_k}} \quad \text{and} \quad \frac{\sum_i \partial m_k^*(\hat{\boldsymbol{\theta}})/\partial \theta_i (\theta_{0,i} - \hat{\theta}_i)}{\sqrt{np_k}}.$$

Also, $(\tilde{\boldsymbol{\theta}} - \hat{\boldsymbol{\theta}})$ is asymptotically normal with mean $\mathbf{0}$ and covariance matrix equal to the negative inverse of the information matrix [Chen (1985)]. So, too, is $(\hat{\boldsymbol{\theta}} - \boldsymbol{\theta}_0)$, and the two quantities are asymptotically independent [e.g., Olver (1974) and Cox and Hinkley (1974)].

Following Chernoff and Lehmann (1954), define $\boldsymbol{\epsilon}$ to be a $K \times 1$ vector with components

$$\epsilon_k = \frac{m_k(\boldsymbol{\theta}_0) - np_k}{\sqrt{np_k}},$$

and let $\hat{\boldsymbol{\nu}}$ be the vector with components

$$\hat{\nu}_k = \sqrt{n}(\hat{p}_k - p_k)/\sqrt{p_k}.$$

It follows from their results that

$$(16) \quad \hat{\boldsymbol{\nu}} = \mathbf{D}(\tilde{\mathbf{J}} + \mathbf{J}^*)^{-1}(\mathbf{D}'\boldsymbol{\epsilon} + \sqrt{n}\mathbf{A}^*) + o_p(1),$$

where $\mathbf{J}^*$ is the matrix whose $(i,j)$th component is

$$\mathbf{E}\left[\frac{\partial \log g(y|\mathbf{z}, \boldsymbol{\theta})}{\partial \theta_i} \frac{\partial \log g(y|\mathbf{z}, \boldsymbol{\theta})}{\partial \theta_j}\right],$$

$g(y|\mathbf{z}, \boldsymbol{\theta})$ is the conditional distribution of $y$ given $\mathbf{z}$ and $\boldsymbol{\theta}$, $\tilde{\mathbf{J}} \equiv \mathbf{D}'\mathbf{D}$ is the matrix with elements

$$\sum_{k=1}^{K} \frac{1}{p_k} \frac{\partial p_k}{\partial \theta_a} \frac{\partial p_k}{\partial \theta_b},$$

and $\mathbf{A}^*$ is the vector whose $a$th component is

$$\frac{1}{n} \sum_{j=1}^{n} \frac{\partial \log g(y|z_j, \theta)}{\partial \theta_a}.$$

From the second corollary, the right-hand side of (16) also describes the large sample distribution of $(m_k(\boldsymbol{\theta}_0) - m_k(\hat{\boldsymbol{\theta}}))/\sqrt{np_k}$.

Taking $\boldsymbol{\eta} = \sqrt{n}\mathbf{A}^*$ and invoking the central limit theorem, Chernoff and Lehmann note that the asymptotic distribution of $(\boldsymbol{\epsilon}, \boldsymbol{\eta})$ is

$$N\left[0, \begin{pmatrix} \mathbf{I} - \mathbf{qq}' & 0 \\ 0 & \mathbf{J}^* \end{pmatrix}\right], \tag{17}$$

where $\mathbf{q}$ is the vector with components $\sqrt{p_k}$. Letting $\boldsymbol{\varepsilon}$ denote a variable having the same distribution as $\boldsymbol{\epsilon}$, and $\boldsymbol{\tau}$ denote a variable having the same distribution as $\boldsymbol{\eta}$, with all four variables distributed independently, it follows that $R^B$ has the asymptotic distribution of

$$(\mathbf{T}\boldsymbol{\varepsilon} + \mathbf{S}\boldsymbol{\tau} - \mathbf{T}\boldsymbol{\epsilon} - \mathbf{S}\boldsymbol{\eta} + \boldsymbol{\epsilon})'(\mathbf{T}\boldsymbol{\varepsilon} + \mathbf{S}\boldsymbol{\tau} - \mathbf{T}\boldsymbol{\epsilon} - \mathbf{S}\boldsymbol{\eta} + \boldsymbol{\epsilon}),$$

where $\mathbf{S} = \mathbf{D}(\tilde{\mathbf{J}} + \mathbf{J}^*)^{-1}$ and $\mathbf{T} = \mathbf{SD}'$. Noting that $\mathbf{D}'\mathbf{q} = \mathbf{0}$, the asymptotic distribution of $(\mathbf{T}\boldsymbol{\varepsilon} + \mathbf{S}\boldsymbol{\tau} - \mathbf{T}\boldsymbol{\epsilon} - \mathbf{S}\boldsymbol{\eta} + \boldsymbol{\epsilon})'$ is $N(\mathbf{0}, \mathbf{I} - \mathbf{qq}')$. The result follows. □

PROOF OF COROLLARY 1. Because the proof of this corollary is similar to the proof of Theorem 1, only an outline is presented here.

To begin, note that Lemma A.1 and Corollary A.2 extend to this setting if $m_k^*(\boldsymbol{\theta})$ is redefined as

$$m_k^*(\boldsymbol{\theta}) = \sum_{j=1}^{n} \mathbf{E}[\text{Ind}(y_j \in [F_j^{-1}(a_{k-1}|\boldsymbol{\theta}), F_j^{-1}(a_k|\boldsymbol{\theta})])].$$

Next, Lemma A.3 applies if (12) is modified so that

$$\begin{aligned}\hat{p}_{j,k} &= F_j[F_j^{-1}(a_k|\boldsymbol{\theta}_0)|\hat{\boldsymbol{\theta}}] - F_j[F_j^{-1}(a_{k-1}|\boldsymbol{\theta}_0)|\hat{\boldsymbol{\theta}}] \\ &= \int_{F_j^{-1}(a_{k-1}|\boldsymbol{\theta}_0)}^{F_j^{-1}(a_k|\boldsymbol{\theta}_0)} f_j(y|\hat{\boldsymbol{\theta}})\,dy,\end{aligned} \tag{18}$$

where $p_{j,k}$ and related estimates refer to the probability that the $j$th observation falls into the $k$th bin. Then

$$\hat{p}_{j,k} - p_{j,k} = (z_{j,k}^*(\boldsymbol{\theta}_0) - z_{j,k}^*(\hat{\boldsymbol{\theta}})) + O_p\left(\frac{1}{n}\right), \tag{19}$$

where

$$z_{j,k}^*(\boldsymbol{\theta}) = \mathbf{E}[\text{Ind}(y_j \in [F_j^{-1}(a_{k-1}|\boldsymbol{\theta}), F_j^{-1}(a_k|\boldsymbol{\theta})])].$$

Corollary A.4 generalizes to

$$\frac{1}{\sqrt{n}}\sum_{j=1}^{n}(\hat{p}_{j,k} - p_{j,k}) = \frac{1}{\sqrt{n}}\sum_{j=1}^{n}(z_{k,j}(\boldsymbol{\theta}) - z_{k,j}(\hat{\boldsymbol{\theta}})) + O_p\left(\frac{1}{\sqrt{n}}\right).$$

Extending Chernoff and Lehmann's (1954) result to the case of nonidentically distributed random variables requires the following modifications of the definitions of variables used in the i.i.d. case. Let

$$\boldsymbol{\epsilon}_j = \left(\frac{z_{j,1} - p_{j,1}}{\sqrt{np_{j,1}}}, \ldots, \frac{z_{j,K} - p_{j,K}}{\sqrt{np_{j,K}}}\right)', \qquad \boldsymbol{\epsilon} = (\boldsymbol{\epsilon}_1', \ldots, \boldsymbol{\epsilon}_n')',$$

$$\tilde{\mathbf{J}} = \left\|\sum_{\alpha=1}^{n}\sum_{r=1}^{K} \frac{1}{p_{\alpha,r}} \frac{\partial p_{\alpha,r}}{\partial \theta_i} \frac{\partial p_{\alpha,r}}{\partial \theta_j}\right\|,$$

$$\mathbf{D} = \begin{pmatrix} \frac{1}{\sqrt{p_{1,1}}}\frac{\partial p_{1,1}}{\partial \theta_1} & \cdots & \frac{1}{\sqrt{p_{1,1}}}\frac{\partial p_{1,1}}{\partial \theta_s} \\ \vdots & \vdots & \\ \frac{1}{\sqrt{p_{1,K}}}\frac{\partial p_{1,K}}{\partial \theta_1} & \cdots & \frac{1}{\sqrt{p_{1,K}}}\frac{\partial p_{1,K}}{\partial \theta_s} \\ \frac{1}{\sqrt{p_{2,1}}}\frac{\partial p_{2,1}}{\partial \theta_1} & \cdots & \frac{1}{\sqrt{p_{2,1}}}\frac{\partial p_{2,1}}{\partial \theta_s} \\ \vdots & \vdots & \\ \frac{1}{\sqrt{p_{n,K}}}\frac{\partial p_{n,K}}{\partial \theta_1} & \cdots & \frac{1}{\sqrt{p_{n,K}}}\frac{\partial p_{n,K}}{\partial \theta_s} \end{pmatrix},$$

$$\mathbf{P} = \Big(\underbrace{\mathbf{I}_k | \ldots | \mathbf{I}_k}_{n \text{ times}}\Big),$$

$$\mathbf{J}^* = \left|\mathbf{E}\left[\left(\sum_{\alpha=1}^n \frac{\partial \log g_\alpha(y|z,\theta)}{\partial \theta_i}\right)\cdot\left(\sum_{\beta=1}^n \frac{\partial \log g_\beta(y|z,\theta)}{\partial \theta_j}\right)\right]\right|,$$

$$A_i^* = \frac{1}{n}\sum_{j=1}^n \frac{\partial \log g_j(y|z,\theta)}{\partial \theta_i}, \qquad \hat{\nu}_{j,r} = \frac{\hat{p}_{j,r} - p_{j,r}}{\sqrt{np_{j,r}}}.$$

Then

$$\hat{\boldsymbol{\nu}} = \mathbf{D}(\tilde{\mathbf{J}} + \mathbf{J}^*)^{-1}(\mathbf{D}'\boldsymbol{\epsilon} + \sqrt{n}\mathbf{A}^*) + o_p(1).$$

The covariance matrix of $\boldsymbol{\epsilon}$ may be written

$$\frac{1}{n}\mathbf{I}_{n\times K} - \frac{1}{n}\begin{pmatrix} \mathbf{q_1 q_1}' & \mathbf{0} & \cdots & \mathbf{0} \\ \vdots & \vdots & \vdots & \vdots \\ \mathbf{0} & \cdots & \mathbf{0} & \mathbf{q_n q_n}' \end{pmatrix},$$

where $\mathbf{q}_i$ is the vector whose $j$th component is $\sqrt{p_{i,j}}$. Denote the rightmost matrix in this equation by $\mathbf{Q}$. Similarly, define $\boldsymbol{\eta} = \sqrt{n}\mathbf{A}^*$. Then the asymptotic distribution of $\boldsymbol{\eta}$ has mean $\mathbf{0}$ and covariance matrix equal to $\mathbf{J}^*/n$, and is independent of $\boldsymbol{\epsilon}$.

Letting $\hat{\mathbf{r}}$ denote the vector with components $(z_{k,j}(\boldsymbol{\theta}) - z_{k,j}(\hat{\boldsymbol{\theta}}))/(\sqrt{np_{j,k}})$, it follows from the generalization of Corollary A.4 that the distribution of $\mathbf{P}\hat{\mathbf{r}}$ is asymptotically the same as that of $\mathbf{P}\hat{\boldsymbol{\nu}}$. Letting $\tilde{\mathbf{r}}$ denote the vector with components $(z_{k,j}(\tilde{\boldsymbol{\theta}}) - z_{k,j}(\hat{\boldsymbol{\theta}}))/(\sqrt{np_{j,k}})$, then $\mathbf{P}\tilde{\mathbf{r}}$ and $\mathbf{P}\hat{\mathbf{r}}$ are, for large $n$, independent and identically distributed. Noting that

$$R^B = (\boldsymbol{\epsilon} - \hat{\mathbf{r}} + \tilde{\mathbf{r}})'\mathbf{P}'\mathbf{P}(\boldsymbol{\epsilon} - \hat{\mathbf{r}} + \tilde{\mathbf{r}})$$

and that $\mathbf{D}'\mathbf{Q} = \mathbf{0}$, some algebra and application of the central limit theorem yields the desired result. $\square$

PROOF OF COROLLARY 2. Expanding the components of $R^B(\tilde{\boldsymbol{\theta}})$ yields

$$(20) \quad \frac{m_k - np_k(\tilde{\boldsymbol{\theta}})}{\sqrt{n}} = \frac{m_k - np_k(\boldsymbol{\theta}_0)}{\sqrt{n}} - \frac{n(p_k(\hat{\boldsymbol{\theta}}) - p_k(\boldsymbol{\theta}_0))}{\sqrt{n}} - \frac{n(p_k(\tilde{\boldsymbol{\theta}}) - p_k(\hat{\boldsymbol{\theta}}))}{\sqrt{n}}.$$

Asymptotically, Taylor series expansions show that the second term on the right-hand side of this equation has the distribution of $\mathbf{T}\boldsymbol{\epsilon} + \mathbf{S}\boldsymbol{\eta}$ described in the proof of Theorem 1, while the third term has the distribution of $\mathbf{T}\boldsymbol{\varepsilon} + \mathbf{S}\boldsymbol{\tau}$. The result follows using methodology in the proof of Theorem 1. $\square$

**Acknowledgments.** The author is grateful to the Editor, an Associate Editor and two referees for helpful comments and suggestions.

DEPARTMENT OF BIOSTATISTICS
AND APPLIED MATHEMATICS
UNIVERSITY OF TEXAS
M. D. ANDERSON CANCER CENTER
1515 HOLCOMBE BLVD., UNIT 447
HOUSTON, TEXAS 77030-4009
USA
E-MAIL: vejohnson@mdanderson.org